\def\[{[\![}
\def\]{]\!]}
\def\Z{{\bf Z}}
\def\C{{\bf C}}
\def\Ch{{\C[[h]]}}
\def\l{\ldots}
\def\n{\noindent}
\def\h{\hat{h}}
\def\e{\hat{e}}
\def\f{\hat{f}}
\def\n{\noindent} 
\def\t{\theta}
\def\s{\smallskip}
\def\b{\bigskip}
\def\a{\hat{a}}
\def\k{{\bar k}}
\def\L{{\bar L}}
\def\q{{\bar q}}

\def\d{\delta}
\def\v{\varepsilon}

\baselineskip=16pt
\font\twelvebf=cmbx12

\noindent
{\twelvebf A  description of the quantum
     superalgebra U$_{\bf q}$[sl(n+1$|$m)] via \hfill\break
     creation and annihilation  generators}

\leftskip 36pt

\vskip 32pt
\noindent
T.D. Palev$^*$ and N.I. Stoilova\footnote*{Permanent address: 
Institute for Nuclear Research and
Nuclear Energy, 1784 Sofia, Bulgaria; E-mail:
tpalev@inrne.bas.bg, stoilova@inrne.bas.bg}

\noindent
Abdus Salam International Centre for Theoretical Physics, 34100 Trieste, Italy

\vskip 1cm

\n{\bf Abstract.} A description of the quantum superalgebra
$U_q[sl(n+1|m)]$ and in particular of the special linear
superalgebra $sl(n+1|m)$ via creation and annihilation generators
(CAGs) is given. It provides an alternative to the canonical
description of $U_q[sl(n+1|m)]$ in terms of Chevalley generators.
A conjecture that the Fock representations of the CAGs provide
microscopic realizations of exclusion statistics is formulated.

\vskip 1cm 

\leftskip 0pt

\n
{\bf 1. Introduction}

\bigskip\n
The description of the quantized simple (universal enveloping) Lie algebras 
[1, 2] and the  basic Lie superalgebras [3-7] is usually carried out in terms 
of their Chevalley generators ($e_i, f_i, h_i, \; i=1,\ldots,n$, for an
algebra of rank $n$). Recently it has been pointed out
that the quantum (super)algebras $U_q[osp(1|2n)]$ [8-10],
$U_q[so(2n+1)]$ [11], more generally 
$U_q[osp(2r+1|2m)],\; r+m=n$ [12],
and also $U_q[sl(n+1)]$ [13] can be defined  via alternative sets 
of generators $a_i^\pm,\; H_i,\;\;i=1,\ldots,n$, referred to as 
(deformed) creation and annihilation generators (CAGs) or
creation and annihilation operators. 

The concept of creation and annihilation generators
of a simple Lie (super)algebra  was introduced in [14].
Let ${\cal A}$ be such an algebra with a
supercommutator $\[\;,\;\]$. The root vectors
$a_1^{\xi},\l,a_n^{\xi}$ of  ${\cal A}$ are said to be creation
($\xi=+$) and annihilation ($\xi=-$) generators of  ${\cal A}$, if 
$$
{\cal A}=lin.env.\{a_i^{\xi},\;\[a_j^{\eta},a_k^{\varepsilon}\]\;|\;
i,j,k=1,\l,n;\; \xi, \eta, \varepsilon =\pm \}, \eqno(1)
$$
so that $a_1^+,\l,a_n^+$ (resp. $a_1^-,\l,a_n^-$) are negative
(resp. positive) root vectors of  ${\cal A}$.

The justification for such terminology stems 
from the observation that the creation and the annihilation
generators of the orthosymplectic Lie superalgebra 
(LS) $osp(2r+1|2m)$ have a direct physical significance: 
$a_1^\pm,\ldots,a_m^\pm$ (resp. $a_{m+1}^\pm,\ldots,a_n^\pm$)
are para-Bose (resp. para-Fermi) operators [15], namely operators which
generalize  the statistics of the tensor (resp. spinor) 
fields in quantum field theory [16]. The LS 
$osp(2r+1|2m)$ is an algebra from the class $B$ in the classification
of Kac [17]. Therefore the paraquantizations (and hence the canonical
Bose and Fermi quantization) could be called $B$-quantizations
(or, more precisely, representations of a $B$-quantization).

A conjecture, stated in [18], assumes that to each class $A$, $B$,
$C$ and $D$ of basic LSs [17] there corresponds a quantum statistics, 
so that its CAGs 
can be interpreted as creation and annihilation operators of real
particles in the corresponding Fock space(s). This conjecture
holds for the classes $A$, $B$, $C$ and $D$ of simple Lie algebras
[19]. It was studied in more details for
the Lie algebras $sl(n+1)$ ($A-$statistics) [20]
and for the LSs $sl(1|m)$ ($A-$superstatistics) [14, 21]. As
an illustration we mention that the
Wigner quantum systems (WQSs), introduced in [22],
are based on the $A-$superstatistics. 
These systems, which attracted some attention
from different points of view [23-25], possess quite
unconventional physical properties. For example, 
the $(n+1)-$particle WQS, based on the LS
$sl(1|3n)$ [26], exhibits a quark like structure: the composite
system occupies a small volume around the centre of mass and
within it the geometry is noncommutative. The underlying
statistics is a Haldane exclusion statistics [27], a subject
of considerable interest in condensed matter physics.

We are not going to discuss further the properties of the
superstatistics (for more details along
this line see [28, 26]  and the references therein). 
We mentioned this point here only in order to indicate that 
the alternative description of $sl(n+1|m)$ and $U_q[sl(n+1|m)]$
will be carried out in terms of (deformed) creation and annihilation 
generators, which, contrary to the Chevalley generators, could be of 
direct physical relevance too.

\bigskip

Throughout the paper we use the notation:

\smallskip
LS, LS's - Lie superalgebra, Lie superalgebras;

CAGs - creation and annihilation generators;

lin.env. - linear envelope;

$\Z$ - all integers;

$\Z_+$ - all non-negative integers;

$\Z_2=\{\bar{0},\bar{1}\}$ - the ring of all integers modulo 2;

$\C$ - all complex numbers;

$[p;q]=\{p,p+1,p+2,\l,q-1,q\}$, for $p\le q\in \Z $;\hfill (2)

\smallskip
$
\t_i=\cases {{\bar 0}, & if $\; i=0,1,2, \ldots , n$,\cr 
               {\bar 1}, & if $\; i=n+1,n+2,\ldots ,n+m$,\cr }; \quad
\t_{ij}=\t_i+\t_j; \hfill (3)
$

\smallskip
$
[a,b]=ab-ba,\;\; \{a,b\}=ab+ba,
\;\;\[a,b\]=ab-(-1)^{deg(a)deg(b)}ba; \hfill (4)
$

$
[a,b]_x=ab-xba,\;\; \{a,b\}_x=ab+xba,
\;\;\[a,b\]_x=ab-(-1)^{deg(a)deg(b)}xba. \hfill (5)
$

\smallskip
\bigskip\n
{\bf 2. The Lie superalgebra sl(n+1$|$m)}

\bigskip\n
Here we give an alternative  definition of the special linear 
Lie superalgebra $sl(n+1|m)$ in terms of creation and
annihilation generators $a_1^\pm,a_2^\pm,\l,a_{n+m}^\pm$. We outline
the relations between the CAGs and the Chevalley generators.

To begin with we recall that the universal enveloping algebra  
$U[gl(n+1|m)]$ of the general linear LS
$gl(n+1|m)$ is a $\Z_2-$graded associative unital superalgebra 
generated by $(n+m+1)^2\;$ $\Z_2-$graded 
indeterminates $\{e_{ij}|i,j\in [0;n+m]\}$, $deg(e_{ij})=\t_{ij}$,
subject to the relations
$$
\[e_{ij},e_{kl}\]=\d_{jk}e_{il}-(-1)^{\t_{ij}\t_{kl}}\d_{il}e_{kj}\quad
i,j,k,l\in [0;n+m]. \eqno(6)
$$
The LS $gl(n+1|m)$ is a subalgebra of $U[gl(n+1|m)]$, considered
as a Lie superalgebra, with
generators  $\{e_{ij}|i,j\in [0;n+m]\}$ and supercommutation
relations (6). The LS $sl(n+1|m)$ is a subalgebra of $gl(n+1|m)$:
$$
sl(n+1|m)=lin.env.\{e_{ij}, (-1)^{\t_k}e_{kk}-(-1)^{\t_l}e_{ll}|
i\ne j;\;  i,j,k,l\in [0;n+m]\}. \eqno(7)
$$
The generators $e_{00},e_{11},\l,e_{n+m,n+m}$ constitute a basis in 
the Cartan subalgebra of \hfill\break
 $gl(n+1|m)$. Denote by $\v_0,\v_1,\l,
\v_{n+m}$ the dual basis, $\v_i(e_{jj})=\d_{ij}$. The root vectors of
both $gl(n+1|m)$ and $sl(n+1|m)$ are $e_{ij}, \; i\ne j,\;
i,j\in [0;n+m]$. The root corresponding to $e_{ij}$
is $\v_i-\v_j$. With respect to the natural order of the basis in the 
Cartan subalgebra $e_{ij}$ is a positive (resp. a negative) root 
vector if $i<j$ (resp. $i>j$).

The above description of $sl(n+1|m)$ is simple, but it is not
appropriate for quantum deformations. A more ``economic'' definition
is given in terms of the Chevalley generators
$$
\h_i=e_{i-1,i-1}-(-1)^{\t_{i-1,i}}e_{ii},\quad \e_i=e_{i-1,i}, \quad
\f_i=e_{i,i-1},\quad i\in [1;n+m] \eqno(8) 
$$
and the $(n+m)\times (n+m)$ Cartan matrix $\{\alpha_{ij}\}$ with
entries
$$
\alpha_{ij}=(1+(-1)^{\t_{i-1,i}})\delta_{ij}-
(-1)^{\t_{i-1,i}}\delta_{i,j-1}-\delta_{i-1,j},\quad 
i,j\in [1;n+m]. \eqno(9) 
$$
We are working with a nonsymmetric Cartan matrix [17]. 
For instance the Cartan matrix (9), corresponding to $n+1=3, \;\;m=5$
is $ 7 \times 7$ dimensional matrix:
$$
(\alpha_{ij})=\pmatrix{
    2&-1& 0& 0& 0& 0& 0\cr
   -1& 2&-1& 0& 0& 0& 0\cr
    0&-1& 0& 1& 0& 0& 0\cr
    0& 0&-1& 2& -1& 0& 0\cr
    0& 0& 0& -1&2& -1& 0\cr
    0& 0& 0& 0& -1&2& -1\cr
    0& 0& 0& 0& 0& -1&2\cr
    }. \eqno(10)
$$
$U[sl(n+1|m)]$ is an associative unital algebra of the
Chevalley generators, subject to the Cartan-Kac relations
$$
\eqalign{
& [\h_i,\h_j]=0,\quad
[\h_i,\e_j]=\alpha_{ij}\e_j,\quad
 [\h_i,\f_j]=-\alpha_{ij}\f_j,\quad
\[\e_i,\f_j\]=\delta _{ij}\h_i,  
}\eqno(11)
$$
and the Serre relations 
$$
\eqalignno{
& [\e_i, \e_j]=0, \quad  [\f_i, \f_j]=0,\quad
if \; |i-j|\neq 1;& (12a) \cr
& \e^2_{n+1}=0,\quad \f^2_{n+1}=0;  & (12b)\cr
& [\e_i, [\e_i, \e_{i+1}]]=0,\quad [\f_i, [\f_i, \f_{i+1}]]=0,
\quad i\ne n+m;& (12c) \cr
& [\e_{i+1}, [\e_{i+1}, \e_{i}]]=0,\quad [\f_{i+1}, [\f_{i+1}, \f_{i}]]=0,
\quad i\ne n+m;& (12d) \cr
& \{[\e_{n+1}, \e_n],[\e_{n+1}, \e_{n+2}]\}=0,\quad
 \{[\f_{n+1}, \f_n],[\f_{n+1}, \f_{n+2}]\}=0. & (12e) \cr
}
$$
The so-called additional Serre relations ($12e$) [29, 30, 31]
can be written also in the form
$$
\{\e_{n+1}, [[\e_n,\e_{n+1}], \e_{n+2}]\}=0,\quad
\{\f_{n+1}, [[\f_n,\f_{n+1}], \f_{n+2}]\}=0. \eqno(12f)
$$

The grading on $U[sl(n+1|m)]$ is induced from the requirement
that the only odd generators are $\e_{n+1}$ and $\f_{n+1}$,
namely
$$
deg(\h_i)={\hat 0},\quad deg(\e_i)=deg(\f_i)=\t_{i-1,i}. \eqno(13)
$$
The LS $sl(n+1|m)$ is a subalgebra of  $U[sl(n+1|m)]$, generated 
by the Chevalley generators in a %%%sence
sense of a Lie superalgebra.
It is a linear span of the Chevalley generators (8)
and all root vectors
$$
\eqalign{
& e_{ij}=[[[\l[[\e_{i+1},\e_{i+2}],\e_{i+3}],\l],\e_{j-1}],\e_j],\cr
& e_{ji}=[\f_j,[\f_{j-1},[\l,[\f_{i+2},\f_{i+1}]\l,]]],
\quad i+1<j; \; i,j\in [0;n+m].\cr
}\eqno(14)
$$

Consider the following root vectors from $sl(n+1|m)$:
$$
\a_i^+=e_{i0}, \quad  \a_i^-=e_{0i}, \;\;i\in [1;n+m], \eqno(15)
$$ 
or, equivalently from (14)
$$
\eqalignno{
& \a_1^-=\e_1,\quad 
\a_i^-=[[[\l[[\e_1,\e_2],\e_3],\l],\e_{i-1}],\e_i]=[\a_{i-1}^-, e_i],
\quad i\in [2;n+m],   & (16a)\cr
&\a_1^+=\f_1, \quad 
\a_i^+=[\f_i,[\f_{i-1},[\l ,[\f_3,[\f_2,\f_1]]\l]]]=[f_i,\a_{i-1}^+].
\quad i\in [2;n+m].   & (16b)\cr
}
$$
The root of $a_i^-$ (resp. of $a_i^+$) is $\v_0-\v_i$ 
(resp. $\v_i-\v_0$). Therefore (with respect to the natural
order of the basis $\v_0,\v_1,\l,\v_{n+m}$) 
$a_1^-,\l,a_{n+m}^-$ are positive root  vectors,
whereas  $a_1^+,\l,a_{n+m}^+$
are negative root vectors. Moreover,  Eq. (1) with 
${\cal A}=sl(n+1|m)$ holds. Hence, the generators (15)
are creation and annihilation generators of $sl(n+1|m)$.
These generators satisfy the following triple relations:
$$
\eqalignno{
& \[\a_i^\xi ,\a _j^\xi \]=0, \quad \xi=\pm, \quad i,j=1,2,\l ,n+m,
& (17a) \cr
%&\cr
& \[ \[ \a_i^+ ,\a _j^- \], \a _k^+\]=\delta_{jk}\a _i^++
(-1)^{\t_i}\delta_{ij}\a _k^+,\quad i,j,k=1,2,\l ,n+m,
& (17b)\cr
%&\cr
& \[ \[ \a_i^+ ,\a _j^- \], \a _k^-\]=-(-1)^{\t_{ij}\t_k}\delta_{ik}
\a _j^- -
(-1)^{\t_i}\delta_{ij}\a _k^-,\quad i,j,k\in [1;n+m].
& (17c)\cr
}
$$
The CAGs (15) together with (17) define completely
$sl(n+1|m)$. The relations (17) are however (similar as Eqs. (6)) not
convenient for quantization. It turns out, and this is a new
result, that one can take only a part of the relations (17), so that they
still define completely $sl(n+1|m)$ and, as we shall see, are appropriate 
for Hopf algebra deformations.

\smallskip\n
{\it Proposition 1.} $U[sl(n+1|m)]$ is an associative unital 
superalgebra with generators 
$\a_i^{\pm}, \;\; i\in [1;n+m]$ and relations:
$$
\eqalignno{
& \[ \a_1^\xi ,\a _2^\xi \]=0,  \quad 
\[ a_1^\xi , a_1^\xi\]=0, \quad \xi=\pm, & (18a) \cr
%&\cr
& \[ \[ \a_i^+ ,\a _j^- \], \a _k^+\]=\delta_{jk}\a _i^++
(-1)^{\t_i}\delta_{ij}\a _k^+,\quad
|i-j|\leq 1, \;\;  i,j,k\in [1;n+m], & (18b) \cr
%&\cr
& \[ \[ \a_i^+ ,\a _j^- \], \a _k^-\]=-(-1)^{\t_{ij}\t_k}\delta_{ik}
\a _j^- -
(-1)^{\t_i}\delta_{ij}\a _k^-,\;; |i-j|\leq 1, \;
  i,j,k\in [1;n+m] & (18c) \cr
}
$$
The $\Z_2-$grading in $U[sl(n+1|m)]$ is induced from
$$
deg(\a_i^\pm)=\t_i. \quad 
\eqno(19)
$$

The proof follows from the expressions of the Chevalley generators 
(8) via the CAGs:
$$
\eqalignno{
& \h_1=\[\a_1^-,\a_1^+\],\quad 
  \h_i=(-1)^{\t_{i-1}}(\[\a_i^-,\a_i^+\]-\[\a_{i-1}^-,\a_{i-1}^+\]),
\quad i\in [2;n+m],& 
(20a)\cr
& \e_1=\a_1^-,\quad \f_1=\a_1^+, \quad
 \e_i=\[\a_{i-1}^+,\a_i^-\],\quad \f_i=\[\a_{i}^+,\a_{i-1}^-\].\quad
 i\in [2;n+m]. & (20b)  \cr
}
$$
We skip the proof of Eqs. (20), since we will give a detailed proof
in the quantum case (see the {\it Theorem}). Only from (18) one 
derives also the larger set of relation (17).

\smallskip
\bigskip\n
{\bf 3. Description of U$_{\bf q}$[sl(n+1$|$m)] via deformed CAGs}

\bigskip\n
In this section we define the quantum superalgebra $U_q[sl(n+1|m)]$
in terms of deformed creation and annihilation generators
$a_i^\pm,H_i,\;\; i=1,2,\l,n+m$. The CAGs are
elements from the so-called Cartan-Weyl basis of $U_q[sl(n+1|m)]$.
A general procedure to construct such a basis
was given in [7] (see also [29]). We follow this procedure and
identify the deformed $a_1^\pm,\l,a_{n+m}^\pm$ generators 
with those elements of the Cartan-Weyl 
basis, which  reduce to the nondeformed CAGs (16) in the limit 
$q \rightarrow 1$. 

First we introduce $U_q[sl(n+1|m)]$ by means of its classical 
definition in terms of the Cartan matrix (9) and the Chevalley
generators.  Let $\Ch$ be the complex algebra of formal power series 
in the indeterminate $h$, $q=e^h\in \Ch$.
$U_q[sl(n+1|m)]$ is a Hopf algebra, which
is a topologically free $\Ch$ module (complete in the $h-$adic
topology), with (Chevalley) generators 
$\{h_i, e_i,f_i\}_{i\in [1;n+m]}$ subject to the
Cartan-Kac relations (${\bar q}=q^{-1}$)
$$
\eqalignno{
& [h_i,h_j]=0,& (21a)\cr
& [h_i,e_j]=\alpha_{ij}e_j,\quad [h_i,f_j]=-\alpha_{ij}f_j,& (21b)\cr
& \[e_i,f_j\]=\delta _{ij}{k_i-\bar{k}_i\over{q-\bar{q}}},\quad
  k_i=q^{h_i},\;k_i^{-1}\equiv\k_i=q^{-h_i}    ,&(21c)\cr 
}
$$
the $e$-Serre relations (see (5))
$$
\eqalignno{
& [e_i,e_j]=0,\; if \; |i-j|\neq 1;\quad  e^2_{n+1}=0; & (22a) \cr
%&&\cr
& [e_i, [e_{i}, e_{i\pm 1}]_{\bar{q}}]_q=
  [e_i, [e_{i}, e_{i\pm 1}]_{q}]_{\bar{q}}=0, \quad i\neq n+1, & (22b)\cr
%&&\cr
& \{ e_{n+1},[[e_n,e_{n+1}]_q, e_{n+2}]_{\bar{q}}\}=
   \{ e_{n+1},[[e_n,e_{n+1}]_{\bar{q}}, e_{n+2}]_{q}\}=0, & (22c)\cr
%&&\cr
}
$$
and the $f-$Serre relations, obtained from the $e$-Serre relations
by replacing everywhere
$e_i$ with $f_i$:
$$
\eqalignno{
& [f_i,f_j]=0,\; if \; |i-j|\neq 1;\quad  f^2_{n+1}=0; & (22d) \cr
%&&\cr
& [f_i, [f_{i}, f_{i\pm 1}]_{\bar{q}}]_q=
  [f_i, [f_{i}, f_{i\pm 1}]_{q}]_{\bar{q}}=0, \quad i\neq n+1, & (22e)\cr
%&&\cr
& \{ f_{n+1},[[f_n,f_{n+1}]_q, f_{n+2}]_{\bar{q}}\}=
   \{ f_{n+1},[[f_n,f_{n+1}]_{\bar{q}}, f_{n+2}]_{q}\}=0. & (22f)\cr
%&&\cr
}
$$

From ($21b$) one derives the following useful relations:
$$
k_ie_j=q^{\alpha_{ij}}e_jk_i, \quad k_if_j=q^{-\alpha_{ij}}f_jk_i,\quad  
\k_ie_j=q^{-\alpha_{ij}}e_j\k_i, \quad \k_if_j=q^{\alpha_{ij}}f_j\k_i.
\eqno(23)
$$

We do not write the other Hopf algebra maps $(\Delta,\; \varepsilon,
S)$ (see [7, 29]), since we will not use them. They are certainly also a
part of the definition.

\smallskip
\n
{\it Remark.} We consider $h$ as an indeterminate. All relations
remain also true, if one replaces $h$ with a number, so that
$q$ is not a root of 1. The latter corresponds to a transition from 
$U_q[sl(n+1|m)$ to the factor algebra $U_q[sl(n+1|m)]/h=number$.

Following [7, 29], introduce a normal order in the 
system of the positive roots
$\Delta_+=\{\v_i-\v_j|i<j\in [0;n+m]\}$ as follows:

\n \centerline{$\v_i-\v_j<\v_k-\v_l\;$ if $\;j<l\;$ or if $\;j=l\;$ 
and $\;i<k. $}

\n Taking into account Eqs. (16), we define the deformed CAGs to be
Cartan-Weyl basis vectors, which are in agreement 
with the above normal order:
$$
\eqalignno{
& a_1^-=e_1,\quad 
a_i^-=[[[\l[[e_1,e_2]_{\bar{q}_1},e_3]_{\bar{q}_2},\l]_{\bar{q}_{i-3}},
e_{i-1}]_{\bar{q}_{i-2}},e_i]_{\bar{q}_{i-1}}=
[a_{i-1}^-, e_i]_{\bar{q}_{i-1}},& (24a) \cr
%& &\cr
&a_1^+=f_1, \quad 
a_i^+=[f_i,[f_{i-1},[\l ,[f_3,[f_2,f_1]_{q_1}]_{q_2}\l]_{q_{i-3}}
]_{q_{i-2}}]_{q_{i-1}}=[f_i,a_{i-1}^+]_{q_{i-1}},&(24b)\cr
& H_1=h_1,\quad H_i=h_1+(-1)^{\t_1}h_2+(-1)^{\t_2}h_3+\l 
+(-1)^{\t_{i-1}}h_i,& (24c)\cr
&&\cr
}
$$
where
$$
q_i=q^{1-2\t_i}=\cases{ q,  & if $i\le n$,\cr
                       \q, & if $i>n$.\cr } \eqno(25)
$$

Note that Eqs. (21)-(23) are invariant with respect to the antilinear 
antiinvolution $(\;)^* $, defined as
$$
(h)^*=-h,\quad  (h_i)^*=h_i,\quad (e_i)^* =f_i, \quad (f_i)^*=e_i,\quad 
(ab)^*=(b)^*(a)^*. \eqno(26)
$$
Therefore
$$
(q)^*={\bar q}, \quad (k_i)^*=\bar{k_i},\quad (a_i^{\pm})^*=a_i^{\mp}, 
\quad (H_i)^*=H_i.   \eqno(27)
$$

The next proposition will be used in several
intermediate computations.

\s\n
{\it Proposition 2.} The relations ($i\ne 1$)
$$
\eqalignno{
& \[ e_i, a_j^-\]_
{q_j^{\delta_{i-1,j}-\delta_{ij}}}=
-q_{i-1}\delta_{i-1,j}a_{i}^-,& (28a) \cr
& \[ f_i, a_j^+\]_
{q_j^{\delta_{i-1,j}-\delta_{ij}}}=
\delta_{i-1,j}a_{i}^+ ,& (28b) \cr
& \[ e_i, a_j^+\]=
\delta_{ij}a_{i-1}^+k_i^{-(-1)^{\t_{i-1}}}, & (28c) \cr
& \[ f_i, a_j^-\]=
-(-1)^{\t_{i-1,i}}\delta_{ij}k_i^{(-1)^{\t_{i-1}}}a_{i-1}^-. & (28d) \cr
}
$$    
follow from (21)-(23) and the definition of the CAGs (24).

\s\n
{\it Proof.}

\n A) Consider first (28a).

\n (i)  The case $j<i-1$. Eq. (28a) is an immediate consequence of
        ($22a$).

\n (ii) The case $j=i-1$ reduces to the definition ($24a$).

\n (iii) The case $j=i$. 

\s       
\n (iii.1) $i=2.$   

\s
\n (iii.1a) If $n=0$, $\[e_2,a_2^-\]_{\q_2}=[e_2,[e_1,e_2]_q]_q
\n            =-q[e_2,[e_2,e_1]_\q]_q=0$, according to ($22b$).

\s
\n (iii.1b)  If $n=1$, 
        $\[ e_2, a_2^-\]_{\bar{q}_2}=\{ e_2, a_2^-\}_{\bar{q}_2}=
        \{e_2,[e_1,e_2]_{\bar{q}}\}_q$

     =$e_2e_1e_2+qe_1e_2^2-\bar{q}e_2^2e_1-\bar{q}qe_2e_1e_2=0$ 
      since, see ($22a$), $e_2^2=0.$ 
 
\s\n (iii.1c)  If $n>1$, $\[ e_2, a_2^-\]_{\bar{q}_2}=
[e_2,[e_1,e_2]_{\bar{q}_1}]_{\bar{q}_2}
=-\bar{q}[e_2, [e_2, e_1]_q]_{\bar{q}}=0$ (see ($22b$).

\s
\n (iii.2) $i>2.$ Using  the identity
$$
{\rm If} \quad \[ a, b\] =0, \;\; 
{\rm then}\;\; \[\[ a,c\]_q , b\]_p=
\[a, \[ c, b\]_p \]_q,\quad p,q\in \Ch,
\eqno(29) 
$$
and the circumstance that $[e_i, a_{i-2}^-]=0$, one obtains
from (24a)

\n
 $a_i^-=[[a_{i-2}^-, e_{i-1}]_{\bar{q}_{i-2}}, e_i]_{\bar{q}_{i-1}}=
 [a_{i-2}^-, [e_{i-1}, e_i]_{\bar{q}_{i-1}}]_{\bar{q}_{i-2}}$.

\s\n (iii.2a) If $i=n+1$,

\n $\[ e_{n+1}, a_{n+1}^-\]_{\bar{q}_{n+1}}=
\{ e_{n+1}, [a_{n-1}^-, [e_n, e_{n+1} ]_{\bar{q}_n}]_{\bar{q}_{n-
1}}\}_{\bar{q}_{n+1}}=             
\{ e_{n+1}, [a_{n-1}^-, [e_n, e_{n+1} ]_{\bar{q}}]_{\bar{q}}\}_q$.

\n Set $ a=e_{n+1}$, $b=a_{n-1}^-$,
$ c=[e_n, e_{n+1}]_{\bar{q}}$; take into account that
$\[ a,b\] =0$ and apply the identity
$$
{\rm If} \;\; \[ a, b\] =0, \quad \[ a, \[ b ,c\]_q\]_p=
(-1)^{\alpha \beta}\[b, \[ a, c\]_p \]_q, \; \alpha=deg(a),\;
\beta=deg(b). \eqno(30) 
$$
\n Then $\[ e_{n+1}, a_{n+1}^-\]_{\bar{q}_{n+1}}=
 [ a_{n-1}^- , z ]_{\bar{q}}$=0, %%%sinse
 since
$z=\{ e_{n+1}, [ e_n ,e_{n+1}]_{\bar{q}} \}_q=0$ 
(follows from $e_{n+1}^2=0$). 

\s
\n (iii.2b) If $i\ne n+1$, then $y= [ e_{i}, [ e_{i} ,e_{i-1}]_{q_{i-1}} ]_
{\bar{q}_{i}}$=0, since in both cases
$i\le n$ or $i>n+1$ it reduces to ($22b$). Therefore, 
$\[ e_i, a_i^-\]_{\bar{q}_{i}}=[ e_i, a_i^-]_{\bar{q}_{i}}=
[ e_i, [a_{i-2}^-, [e_{i-1}, e_{i} ]_{\bar{q}_{i-1}}]_{\bar{q}_{i-
2}}]_{\bar{q}_{i}}$

\n (if $a=e_i,\;\; b=a_{i-2}^-,\;\;
c=[e_{i-1}, e_{i}]_{\bar{q}_{i-1}}$
then $[ a,b] =0$  and  from  (30)) 
 
\n $= [ a_{i-2}^- , [ e_{i}, [ e_{i-1} ,e_{i}]_{\bar{q}_{i-1}} ]_
{\bar{q}_{i}} ]_{\bar{q}_{i-2}}=
-\bar{q}_{i-1}[ a_{i-2}^- ,y]_{\bar{q}_{i-2}}=0$. Hence (28a) holds
for any $i=j>1$.

\s
\n 
(iv) The case $j=i+1.$

\s\n
(iv.1) If $i=2,\;\; n+1\neq 2,$ 
$\[ e_2, a_3^-\]=[e_2,[\;[ e_1,e_2]_{\bar{q}_1}, e_3]_{\bar{q}_2}]=
[e_2,[\;[ e_1,e_2]_{\bar{q}_1}, e_3]_{\bar{q}_1}]$.
\n For $b=e_2,\; a=e_1,\; c=e_3$ use the identity:

If $b$ is even and $ \[ a, c\] =0$,  then
$$
(x+\bar{x})[ b,\[  a, [b,c ]_x \]_x]=\[ a, [ b, [ b, 
c]_x]_{\bar{x}}\]_{x^2} -
\[[ b, [ b, a]_x]_{\bar{x}}, c\]_{x^2}, \quad \bar{x}=x^{-1}.
\eqno(31) 
$$
Then $\[ e_2, a_3^-\]=(q_1+\bar{q}_1)^{-1}\left( \[ e_1, [e_2, [e_2, e_3]_
{\bar{q}_1}]_{q_1}\]_{q_1^{-2}}-\[ [e_2, [e_2, e_1]_
{\bar{q}_1}]_{q_1}, e_3\]_{q_1^{-2}}\right) =0$ according to ($22b$).

\s\n
(iv.2) If $i=2,\;\; n+1=2,\;\; 
\[ e_2, a_3^- \] =\{ e_2, [[e_1, e_2]_{\bar{q}}, e_3]_q\}=0$ 
according to ($22c$).

\b\n (iv.3) For  $i>2,$ set (see (24a)) $a_{i+1}^-=[[[a_{i-2}^-, e_{i-
1}]_{\bar{q}_{i-2}},
e_i]_{\bar{q}_{i-1}}, e_{i+1}]_{\bar{q}_i}.$ Use
that $[a_{i-2}^-,e_i]=[a_{i-2}^-,e_{i+1}]=0$ and apply twice (29):
$a_{i+1}^-= [a_{i-2}^-,[[e_{i-1}, e_i]_{\bar{q}_{i-1}}, 
e_{i+1}
]_{\bar{q}_i}]_{\bar{q}_{i-2}}$.

\s
\n (iv.3a) If $i=n+1$,
$\[ e_i, a_{i+1}^-\] =\{ e_{n+1}, a_{n+2}^-\}=\{ e_{n+1}, [a_{n-1}^-, 
[[e_n, e_{n+1}]_{\bar{q}_n}, e_{n+2}]_{\bar{q}_{n+1}}]_{\bar{q}_{n-
1}}\}$

\n (use that $[e_{n+1},a_{n-1}^-]=0$ and (30))

\n $=[a_{n-1}^-, \{ e_{n+1}, [[ e_n, e_{n+1}]_{\bar{q}_n},
e_{n+2}]_{\bar{q}_{n+1}}\} ]_{\bar{q}_{n-1}}=0$ according to ($22c$)
and (25).

\s\n
(iv.3b) If $i\neq n+1$ $\;\; \[ e_i , a_{i+1}^-\] =[e_i, a_{i+1}^-]=
[e_i, [a_{i-2}^-, [[ e_{i-1}, e_i]_{\bar{q}_{i-1}}, e_{i+1}
]_{\bar{q}_i}]_{\bar{q}_{i-2}}]$

\n ($ [e_i,a_{i-2}^-]=0$, use (30))
 
\n $=[a_{i-2}^-, [e_{i}, [[ e_{i-1}, e_i]_{\bar{q}_{i-1}}, e_{i+1}
]_{\bar{q}_i}]]_{\bar{q}_{i-2}}= 
[a_{i-2}^-, [e_{i}, [[ e_{i-1}, e_i]_{\bar{q}_{i}}, e_{i+1}
]_{\bar{q}_i}]]_{\bar{q}_{i-2}}$.

\n
If $a=e_i$, $b=e_{i-1}$,  $c=e_{i+1}$, then  $[b,c]=0$; 
apply a similar to (31) identity:

If $a$  is even and $\[ b, c\] =0$,  then 
$$
(x+\bar{x})[ a,\[ [ b,a ]_x,c \]_x]=\[ b, [ a, [ a, c]_x]_{\bar{x}}\]_{x^2} -
\[[ a, [ a, b]_x]_{\bar{x}}, c\]_{x^2}. \eqno(32) 
$$
The latter yields
$ \[ e_i , a_{i+1}^-\]$ 

\n $=(\bar{q}_i+q_i)^{-1} [a_{i-2}^-, \left( [e_{i-1},
[e_i, [e_i, e_{i+1}]_{\bar{q}_i}]_{q_i}]_{q_i^{-2}}-
[[e_i,[e_i,e_{i-1}]_{\bar{q}_i}]_{q_i}, e_{i+1}]_{q_i^{-2}}\right) 
]_{\bar{q}_{i-2}}
 =0$ according to ($22b$).

\s\n
(v) The case $j>i+1.$ Then
$a_j^-=[[[\ldots [[a_{i+1}^-, e_{i+2}]_{\bar{q}_{i+1}}, 
e_{i+3}]_{\bar{q}_{i+2}},\ldots ]_{\bar{q}_{j-3}}, e_{j-
1}]_{\bar{q}_{j-2}}
, e_j]_{\bar{q}_{j-1}}$
and since $e_i$ commutes with $e_{i+2}, \;\; e_{i+3}, \;\;\ldots , 
e_j\;$, see ($22a$), and $e_i$ supercommutes with $a_{i+1}^-$, 
see (iv), one concludes that $\[ e_i, a_j^-\]=0.$ 
The unification of (i)-(v) %%%yealds
yields (28a). 

\b\n
B) Applying the antiinvolution (26) on both 
sides of ($28a$) one obtains ($28b$).

\bigskip\n
C) We pass to prove ($28c$).

\bigskip\n
(i) For $i>j$, ($28c$) is an immediate consequence of ($24b$) and 
($21c$).

\bigskip\n
(ii) Let $i=j$. 
$\[ e_i,a_i^+\]=\[ e_i, [f_i, a_{i-1}^+]_{q_{i-1}}\]$

\n(from (i) $ \[e_i, a_{i-1}^+\] =0$, apply  (29))

\n
$=\[\[ e_i, f_i\], a_{i-1}^+\]_{q_{i-1}}
=[ {k_i-\bar{k}_i\over{q-\bar{q}}} , a_{i-1}^+]_{q_{i-1}}=
a_{i-1}^+k_i^{-(-1)^{\t_{i-1}}}$.
\n In the last step we used the relations 
$k_i a_{i-1}^+=qa_{i-1}^+k_i$ and $\k_i a_{i-1}^+=\q a_{i-1}^+\k_i$,
which follow from ($24b$) and ($23$).

\bigskip
\n (iii) Let $j=i+1$. 
$\[e_i, a_{i+1}^+\]= \[e_i,[f_{i+1},[f_i, a_{i-1}^+]_{q_{i-
1}}]_{q_i}\]$

\n (take into account that $[e_i, f_{i+1}]=0$ and apply (30)) 

\n
$ =\[ f_{i+1}, \[ e_i,[f_i, a_{i-1}^+]_{q_{i-1}}\]\]_{q_i}$
(now $[e_i,a_{i-1}^+]=0$, use (29))  

\n $= \[ f_{i+1}, \[\[ e_i,f_i\] , a_{i-1}^+\]_{q_{i-1}}\]_{q_{i}}=
\[ f_{i+1}, \[ {k_i-\bar{k}_i\over{q-\bar{q}}}, a_{i-1}^+\]_{q_{i-
1}}\]_{q_i}=
[ f_{i+1}, a_{i-1}^+k^{-(-1)^{\t_{i-1}}}]_{q_i}$ 

\n Using the identity
$$
[a,bc]_x=[a,b]c+b[a,c]_x \eqno(33)
$$
one has
$\[e_i, a_{i+1}^+\]=[f_{i+1}, a_{i-1}^+]k_i^{-(-1)^{\t_{i-1}}}
+a_{i-1}^+[f_{i+1}, k_i^{-(-1)^{\t_{i-1}}}]_{q_i}=0$, according to 
($28b$), ($23$) and (25).

\bigskip\n
(iv) For $j>i+1\;$
$a_j^+=[f_j,[f_{j-1}, [\ldots ,[ f_{i+3}, [ f_{i+2}, 
a_{i+1}^+]_{q_{i+1}}]_{q_{i+2}}\ldots ]_{q_{j-3}}]_{q_{j-2}}]_{q_{j-1}}$

\n and since $e_i$ supercommutes with $a_{i+1}^+$, see (iii), and 
commutes with $f_j,\; f_{j-1}, \; \ldots , f_{i+2}$, see ($21c$), 
one concludes that $\[ e_i , a_j^+\]=0$.
The unification of (i)-(iv) yields ($28c$).

\s\n
D) Applying the antiinvolution (26) on both 
sides of ($28c$) one obtains ($28d$). 

\b\n This completes the proof.

\bigskip\n
{\it Proposition 3.} The deformed CAG's (23) generate $U_q[sl(n+1|m)].$

\smallskip\n
{\it Proof.} Let
$$
L_i=q^{H_i}, \quad \L_i\equiv L_i^{-1}=q^{-H_i}.\eqno(34)
$$ 
The proof is a consequence of the relations
$$
\eqalignno{
&\[ a_i^-, a_i^+\]={{L_i-\L_i}\over{q-\q}}& (35a)\cr
%&&\cr
&\[ a_i^-, a_{i+1}^+\]=-(-1)^{\t_i}L_if_{i+1} & (35b)\cr
%&&\cr
&\[ a_{i+1}^-, a_{i}^+\]=-(-1)^{\t_i}e_{i+1}\L_i  & (35c)\cr
}
$$    
We prove these equations by induction on $i.$ For $i=1,$ 
$(35a)$ holds. 
Let $(35a)$ be true. Then  from $(28d)$, (30) and $(35a)$ one has

\b\n
$
\[a_i^-, a_{i+1}^+\]=\[a_i^-,[f_{i+1}, a_i^+]_{q_i}\]=\[ f_{i+1},
\[a_i^-, a_i^+\]\]_{q_i}={1\over{q-\q}}[f_{i+1},L_i-\L_i]_{q_i}
$

\b\n
$
={1\over{q-\q}}[f_{i+1}, k_1k_2^{(-1)^{\t_1}}k_3^{(-1)^{\t_2}}
\ldots k_i^{(-1)^{\t_{i-1}}}-\k_1k_2^{-(-1)^{\t_1}}k_3^{-(-
1)^{\t_2}}
\ldots k_i^{-(-1)^{\t_{i-1}}}]_{q_i}.
$

\b\n
Using (25) and repeatedly $(23)$, one end with 

\smallskip\n
$
\[ a_i^-, a_{i+1}^+\]=-(-1)^{\t_i}k_1k_2^{(-1)^{\t_1}}k_3^{(-
1)^{\t_2}}
\ldots k_i^{(-1)^{\t_{i-1}}}f_{i+1},
$
namely with $(35b)$.
Similarly, one proves $(35c)$. Therefore, if $(35a)$ holds, then also 
equations $(35b)$ and $(35c)$ are fulfilled. Assuming this, consider 

\s\n
$\[ a_{i+1}^-, a_{i+1}^+\]=
\[ [a_i^-, e_{i+1}]_{\bar{q}_i}, a_{i+1}^+\].$ 
Then the identity 
$$
\[\[a,b\]_x, c\]=(-1)^{\beta\gamma}\[\[a,c\],b\]_x+
\[a, \[b,c\]\]_x,\quad  \beta=deg(b),\; \gamma=deg(c) \eqno(36)
$$
yields
\b\n
$$
\eqalign{
& \[ a_{i+1}^-, a_{i+1}^+\]
=(-1)^{\t_{i,i+1}} \[\[ a_i^-, a_{i+1}^+\], e_{i+1}\]_{\bar{q}_i}+
\[ a_i^-, \[ e_{i+1}, a_{i+1}^+
\] \]_{\bar{q}_i}\cr
&\cr
& =-(-1)^{\t_{i+1}}[k_1k_2^{(-1)^{\t_1}}k_3^{(-1)^{\t_2}}
\ldots k_i^{(-1)^{\t_{i-1}}}f_{i+1}, e_{i+1}]_{\bar{q}_i}+ 
\[ a_i^-, a_i^+ k_{i+1}^{-(-1)^{\t_i}} \] _{\bar{q}_i}\cr
&\cr
& =-(-1)^{\t_{i+1}} \[ f_{i+1}, e_{i+1}\] k_1k_2^{(-1)^{\t_1}}k_3^{(-
1)^{\t_2}}
\ldots k_i^{(-1)^{\t_{i-1}}}+
\[ a_i^-, a_i^+\] k_{i+1}^{-(-1)^{\t_i}}\cr
&\cr
& ={k_1k_2^{(-1)^{\t_1}}k_3^{(-1)^{\t_2}}
\ldots k_i^{(-1)^{\t_{i-1}}}k_{i+1}^{(-1)^{\t_i}}-
{\bar{k}_1k_2^{-(-1)^{\t_1}}}k_3^{-(-1)^{\t_2}}
\ldots k_i^{-(-1)^{\t_{i-1}}}k_{i+1}^{-(-1)^{\t_i}}\over {q-\bar{q}}}\cr
&\cr
&={{L_{i+1}-\L_{i+1}}\over{q-\q}}.\cr
}
$$

\bigskip
Thus, Eqs (35) hold for any $i$. From (24c) and (35) we have 
$$
\eqalignno{
& e_1=a_1^-,\quad e_{i+1}=-(-1)^{\t_i}\[a_{i+1}^-,a_i^+\]q^{H_i},
\quad i\in [1;n+m-1] & (37a) \cr
& f_1=a_1^+,\quad f_{i+1}=-(-1)^{\t_i}\q^{H_i}
  [a_{i}^-,a_{i+1}^+\], \quad i\in [1;n+m-1] & (37b) \cr
& h_1=H_1,\quad h_i=(-1)^{\t_{i-1}}(H_i-H_{i-1})\quad i\in [2;n+m],
& (37c)
}
$$
which completes the proof.

\b
We proceed to state our main result.

\bigskip\n
{\it Theorem.} $U_q[sl(n+1|m)]$ is an unital associative
algebra, which is topologically free $\Ch$ module, with generators
$\{H_i,\; a_i^\pm\}_{i\in [1;n+m]}$ and relations
$$
\eqalignno{
& [H_i,H_j]=0, & (38a) \cr
%&&\cr
& [H_i,a_j^{\pm}]=\mp(1+(-1)^{\t_i}\delta_{ij})a_j^{\pm}, 
&(38b)\cr
%&&\cr
& \[ a_i^-, a_i^+\]={L_i-\bar{L}_i\over{q-\bar{q}}}, & (38c)\cr
%&&\cr
& \[\[a_i^{\eta}, a_{i+\xi}^{-\eta}\], a_k^{\eta} 
\]_{q^{\xi (1+(-1)^{\t_i}\delta_{ik})}}=
\eta^{\t_k}\delta_{k,i+\xi}L_k^{-\xi\eta}
a_i^{\eta}, \;\;\xi, \eta =\pm \;\;or \;\; \pm 1, & (38d)\cr
%&&\cr
& \[ a_1^\xi , a_2^\xi \]_q =0, \quad \[ a_1^\xi , a_1^\xi \] =0,
\quad \xi=\pm. & (38e) \cr
}
$$
 
\b\n
{\it Proof.} As a first step one has to show that Eqs. (38) hold. 
Most of the results for this part of the proof are already
obtained. 
Eq. ($38a$) is evident. 
Eq. ($38b$) follows from the relation
    $\sum_{p=1}^i\sum_{q=1}^j (-1)^{\t_{p-1}}\alpha_{pq}=
     1+(-1)^{\t_i}\delta_{ij}$, the definitions of $a_i^\pm$ and
     $H_i$ (see (24)) and the relations (21b). From (38b) one
     also derives

$$ 
L_ia_j^{\pm}=q^{\mp (1+(-1)^{\t_i}\delta_{ij}})a_j^{\pm}L_i.\eqno(39)
$$
Eq. ($38c$) is the same as ($35a$). 
The derivation of all triple relations ($38d$) is relatively long, 
but simple. 
It is based  on a case by case considerations. To this end one replaces 
$e_i$ and $f_i$ in (28) with the right hand sides of ($37a,b$). 
The nontrivial part is to put all cases in the compact form ($38d)$.
If $n\ne 0$, $\[ a_1^\xi , a_1^\xi \]=[ a_1^\xi , a_1^\xi ]=0$. The
first relations in $(38e)$ reduce to the triple Serre relations
$(22b,e)$. If $n=0$, Eqs. $(38e)$ hold because $e_1^2=0$ and
$f_1^2=0$.

It remains to prove as a second step that 
any other relation in $U_q[sl(n+1|m)]$ follows from Eqs. (38). 
To this end it suffices to show that all Cartan-Kac relations (21) 
and the Serre relations (22) follow from (38). 

\s\n
A) The Cartan-Kac relations $(21a)$ follow in an %%%eviden
evident way from
   (37c) and (38a).

\s\n
B) Eqs. $(21b)$ are easily derived from (37) and (38b).

\s\n
C) The proof of $(21c)$ is not trivial.

\n
(i) The case $i=j=1$ is evident.

\n
(ii) The case $i=1,\; j>1$:
$
\[ f_j, e_1\] 
=\[ -(-1)^{\t_{j-1}} \bar{L}_{j-1} \[ a_{j-1}^-, a_j^+\],
a_1^-\]$ 

\n
(using (39))

\n$
=-(-1)^{\t_{j-1}}\bar{L}_{j-1}\[\[ a_{j-1}^-, a_j^+\], a_1^-
\]_{q^{(1+(-1)^{\t_{j-1}}\delta_{j-1,1})}}=0 $  according to
$(38d)$.

\n
(iii) In a similar way one shows that $\[ e_i, f_1\]=0$ for 
$i>1.$

\n
(iv) The case $i,j\in [2;n+m]$. From (37)

$
\[ e_i, f_j\]
=(-1)^{\t_{i-1,j-1}}\[ \[ a_i^-, a_{i-1}^+\]L_{i-1},
\bar{L}_{j-1}\[ a_{j-1}^-, a_j^+\]\]$. 

\n
Apply (39):
$$
\eqalign{
 \[ e_i, f_j\]&
=(-1)^{\t_{i-1,j-1}} q^{(-1)^{\t_{i-1}-(-1)^{\t_{j-1}}\delta_{ij}
+(-1)^{\t_{j-1}}\delta_{i,j-1}}}L_{i-1}\bar{L}_{j-1}  \cr
&\times\[\[a_i^-,a_{i-1}^+\], \[ a_{j-1}^-, a_j^+\]\]_{q^{(-1)^{\t_{i-1}}
\delta_{i-1,j}-(-1)^{\t_{j-1}}\delta_{i,j-1}}}\cr
} \eqno(40)
$$
(iv.1) For $i=j$ (40) reduces to

$$
\[ e_i, f_i\] =\[\[ a_i^-,a_{i-1}^+\], \[ a_{i-1}^-, a_i^+\]\].\eqno(41) 
$$
In order to evaluate the r.h.s. of (41) use the following identity
($ \alpha=deg(a),\;\beta=deg(b)$):

If  $x=zs,\; y=zr,\; t=zsr; \quad x,y,z,r,s,t\in \Ch$, then
$$
\[ a,\[ b,c\]_x\]_y=\[\[ a,b\]_z,c\]_t+z(-
1)^{\alpha\beta}\[b,\[a,c\]_r\]_s.\eqno(42)
$$

\n Applying  (42) to the r.h.s. of (41) with
$a=\[ a_{i-1}^+, a_i^-\]$,
$b=a_i^+$, $c=a_{i-1}^-$ and  $x=y=1$,
$z=q$, $r=s=t=\bar{q}$,
one obtains

\n
$
\[ e_i, f_i\] =\[\[\[a_{i-1}^+,a_i^-\],a_i^+\]_q,a_{i-1}^-\]_{\bar{q}}
-q(-1)^{\t_i}\[ a_i^+, \[\[ a_i^-, a_{i-1}^+\], a_{i-1}^-\]_{\bar{q}}
\]_{\bar{q}}\]$

\n (use  ($38d$))

\n
$
=\[ \bar{L}_ia_{i-1}^+, a_{i-1}^-\]_{\bar{q}}-
q(-1)^{\t_{i,i-1}}\[ a_i^+, \bar{L}_{i-1}a_i^-\]_{\bar{q}}
$

\n (use  (39))

\n $
=\bar{L}_i\[ a_{i-1}^+, a_{i-1}^-\]-(-1)^{\t_{i,i-1}}\bar{L}_{i-1}\[ 
a_i^+, a_i^-\]
$

\n (use  $(38c)$  and (34))

\n
$
={(-1)^{\t_{i-1}}\over{q-\bar{q}}}\left( k_i^{(-1)^{\t_{i-1}}}-
k_i^{-(-1)^{\t_{i-1}}}\right). 
$

\n Taking into account that $\t_{i-1}=0,$ for $i\in [1;n+1]$
and that $\t_{i-1}=1$ for $i\in [n+2;n+m]$, one 
ends with
%(see p. $sl(n+1/m)-87)$

$$
\[ e_i, f_i\]={k_i-\bar{k}_i\over{q-\bar{q}}}.\eqno(43)
$$
(iv.2) Let $|i-j|>1.$
Then (40) reduces to 

\n
$
\[ e_i, f_j\]=(-1)^{\t_{i-1,j-1}}q^{(-1)^{\t_{i-1}}}
L_{i-1}\bar{L}_{j-1}\[\[a_i^-,a_{i-1}^+\], \[ a_{j-1}^-, a_j^+\]\]
$

\n
$
=
(-1)^{\t_{ij}}q^{(-1)^{\t_{i-1}}}
L_{i-1}\bar{L}_{j-1}\[\[a_{i-1}^+,a_{i}^-\], \[ a_{j}^+, a_{j-1}^-\]\]
$

\n
(apply (42) with $a=\[ a_{i-1}^+, a_i^-\]$,  $b=a_j^+$,
$c=a_{j-1}^-$, $x=y=1$, $z=q$, $t=s=r=\bar{q}$)

\n
$
=(-1)^{\t_{ij}}q^{(-1)^{\t_{i-1}}}L_{i-1}\bar{L}_{j-1}
\big( \[\[\[ a_{i-1}^+, a_i^-\],a_j^+\]_q, a_{j-1}^-\]_{\bar{q}}
$

\n
$
-q(-1)^{\t_{i,i-1}\t_j+\t_{i,i-1}} \[ a_j^+, \[\[ a_i^-, a_{i-1}^+\], 
a_{j-1}^-\]_{\bar{q}}\]_{\bar{q}}\big)=0
$ from ($38d$).

\s\n(iv.3) For $j=i-1$ $(40)$ reduces to

\n
$
\[ e_i, f_{i-1}\]=(-1)^{\t_{i-1,i-2}}
q^{(-1)^{\t_{i-1}}}L_{i-1}\bar{L}_{i-2}
\[\[ a_{i-1}^+, a_i^-\], \[ a_{i-1}^+, a_{i-2}^-\]\]_{q^{(-1)^{\t_{i-
1}}}}
$

\n
(use (42) with $a=\[ a_{i-1}^+, a_i^-\]$, $b=a_{i-1}^+$,
$c=a_{i-2}^-$, $x=1$, 
$y=q^{(-1)^{\t_{i-1}}}$, 
$r=t=\bar{q}$,
$z=q^{1+(-1)^{\t_{i-1}}}$, 
$s=q^{-(1+(-1)^{\t_{i-1}})}$)

\n
$
=(-1)^{\t_{i-1,i-2}}
q^{(-1)^{\t_{i-1}}}L_{i-1}\bar{L}_{i-2}
( \[\[\[a_{i-1}^+, a_i^-\], a_{i-1}^+\]_{q^{1+(-1)^{\t_{i-1}}}},
a_{i-2}^-\]_{\bar{q}}
$

\n
$
+q^{1+(-1)^{\t_{i-1}}}\[ a_{i-1}^+, \[\[ a_{i-1}^+, a_i^-\], 
a_{i-2}^-\]_{\bar{q}}\]_{q^{-(1+(-1)^{\t_{i-1}})}})=0
$ 
from $(38d)$ 

\s\n
(iv.4) Let $j=i+1$. Then $(40)$ reduces to

\n
$
\[ e_i, f_{i+1}\]= (-1)^{\t_{i-1,i}}q^{(-1)^{\t_{i-1}}+(-1)^{\t_i}}
L_{i-1}\bar{L}_i
\[\[a_i^-, a_{i-1}^+\], \[ a_i^-, a_{i+1}^+\]\]_{q^{-(-1)^{\t_i}}}
$

\n
$
= (-1)^{\t_{i,i+1}}q^{(-1)^{\t_{i-1}}+(-1)^{\t_i}}
L_{i-1}\bar{L}_i
\[\[a_{i-1}^+, a_{i}^-\], \[ a_{i+1}^+, a_{i}^-\]\]_{q^{-(-1)^{\t_i}}}
$

\n (from (42) with $a=\[ a_{i-1}^+, a_i^-\]$, 
$b=a_{i+1}^+$, $c=a_i^-$, $x=1$, $y=q^{-(-1)^{\t_i}}$,  
$z=q$, $s=\bar{q}$, $r=t=q^{-(1+(-1)^{\t_i})}$)

\n
$
= (-1)^{\t_{i,i+1}}q^{(-1)^{\t_{i-1}}+(-1)^{\t_i}}
L_{i-1}\bar{L}_i
( \[\[\[ a_{i-1}^+, a_i^-\], a_{i+1}^+\]_q, a_i^-\]_{q^{-(1+(-
1)^{\t_i})}}
$

\n
$
+q(-1)^{\t_{i-1,i}}\[ a_{i+1}^+, \[\[ a_{i-1}^+, a_i^-\], 
a_i^-\]_{q^{-(1+(-1)^{\t_i})}}\]_{\bar{q}})=0$, according to   
($38d$)).

\b
So far we have shown that all Cartan-Kac relations (21) follow
from (38). It remains to verify the Serre relations (22). 
We consider in some more %%%detais
details the $e-$Serre relations
($22a-c$).

\s\n
D) We pass to prove first  ($22a$), namely that $\[e_i,e_j\]=0$ 
if $|i-j|\ne 1$.
\s\n
(i) The case with $i=1$ and  $j=[3;n+m]$ follows directly from
     (39) and ($38d$).

\s\n
(ii) $i\ne j\in [2;n+m]$. From ($37a$)

\n
$
[e_i, e_j]=(-1)^{\t_{i-1,j-1}}[\[a_i^-,a_{i-1}^+\]
L_{i-1}, \[ a_j^-, a_{j-1}^+\]L_{j-1}]
$

\n
(use ($39$))

\n
$
=(-1)^{\t_{ij}}[\[ a_{i-1}^+, a_i^-\], \[ a_{j-
1}^+,
a_j^-\]]L_{i-1}L_{j-1}
$

\n (apply (42) with $a=\[ a_{i-1}^+, a_i^-\],$ $b=a_{j-1}^+$, 
$c=a_j^-,$ $x=y=1,$ $z=q$, $t=r=s=\bar{q}$)

\n
$
=(-1)^{\t_{ij}}(\[\[\[ a_{i-1}^+, a_i^-\], a_{j-1}^+
\]_q, a_j^-\]_{\bar{q}}+
q(-1)^{\t_{i,i-1}\t_{j-1}}\[ a_{j-1}^+, \[\[ a_{i-1}^+, a_i^-\], a_j^-
\]_{\bar{q}}\]_{\bar{q}})L_{i-1}L_{j-1}=0$, according to ($38d$).

\s\n
(iii) If $i=j\ne n+1$,  $\[e_i,e_i\]=[e_i,e_i]=0$   

\s\n
(iv) Consider $e_{n+1}^2={1\over 2}\{e_{n+1},e_{n+1}\}
     \equiv {1\over 2}\[e_{n+1},e_{n+1}\].$

\n
(iv.1) The case with $n+1=1$ is evident:
$ \{ e_1, e_1\}=\{a_1^-, a_1^-\}=0$, see $(38e)$.

\n
(iv.2) $n+1\neq 1$. Use (37a):
$
e_{n+1}^2\sim \{ e_{n+1}, e_{n+1}\}_{q^2}
=\{ \[ a_{n+1}^-, a_n^+\] L_n, \[ a_{n+1}^-, a_n^+\] L_n\}_{q^2}
$

\n
$
=\bar{q}\[\[ a_{n+1}^-, a_n^+\], \[ a_{n+1}^-
, a_n^+\]\]_{q^2}L_n^2
$

\n
(apply (42) with $a=\[ a_{n+1}^-, a_n^+\]$, $b=a_{n+1}^-$,
$c=a_n^+$, $x=s=z=1$, $y=r=t=q^2$)

\n
$
=\bar{q}(\[\[\[ a_{n+1}^-, a_n^+\], a_{n+1}^-\], a_n^+\]_{q^2}
-\[ a_{n+1}^-, \[\[ a_{n+1}^-, a_n^+\], a_n^+\]_{q^2}\])L_n^2=0$,  
according to ($38d$)). Hence the Serre relations (22a) follow
from (38).

\s\n
E) We prove the triple Serre relation $[e_i, [e_{i}, e_{i+1}]_{\bar{q}}]_q=
  [e_i, [e_{i}, e_{i+1}]_{q}]_{\bar{q}}=0, \quad i\neq n+1.$

\s\n (i) Let $i=1$. Since $n+1\ne 1$,
$\; a_1^-$ is an even generator. Taking this into account,
one easily derives only from ($38$) and (39) that
 
\n $[e_1,e_2]_\q = 
[a_1^-,[a_1^+,a_2^-]L_1]_\q = [a_1^-,[a_1^+,a_2^-]]_q L_1 =a_2^-$.
Therefore, see  ($38e$), \hfill\break
$[e_1,[e_1,e_2]_\q]_q=[a_1^-,a_2^-]_q = 0$.
 
\n\n
(ii) $i\in [2; n].$ From ($37a$) and (39)
$[e_i, e_{i+1}]_{\bar{q}}= [[ a_i^-, a_{i-1}^+] L_{i-1}, 
[ a_{i+1}^-, a_i^+]L_i]_{\bar{q}}\hfill\break 
=[[ a_i^-, a_{i-1}^+], [ a_{i+1}^-, a_i^+]]L_iL_{i-1}$

\n
(apply  (42)  with  $a=[ a_i^-, a_{i-1}^+]$,  $b=a_{i+1}^-$,
$c=a_i^+$, $x=y=1$, $z=\bar{q}$, $r=s=t=q$)

\n $=[[[ a_i^-, a_{i-1}^+], a_{i+1}^-]_{\bar{q}}, a_i^+]_qL_iL_{i-1}
+\bar{q}[ a_{i+1}^-, [[ a_i^-, a_{i-1}^+], a_i^+]_q]_q L_iL_{i-1}$

\n (use  ($38d$) and (39)) 

\n
$ 
=-\bar{q}[ a_{i+1}^-, \bar{L}_ia_{i-1}^+]_qL_iL_{i-1} 
=-[a_{i+1}^-, a_{i-1}^+]L_{i-1}.
$

\s\n
Therefore

\n
$
[e_i, [e_i, e_{i+1}]_{\bar{q}}]_q=[[ a_i^-, a_{i-1}^+] L_{i-1}, [ 
a_{i+1}^-,a_{i-1}^+]L_{i-1}]_q
=\bar{q}[[a_i^-, a_{i-1}^+],[a_{i+1}^-, a_{i-1}^+]]_qL_{i-1}^2
$

\n (from  (42) with $a=[a_i^-, a_{i-1}^+]$, $b=a_{i+1}^-$
, $c=a_{i-1}^+$, $x=1$, $y=s=q$, $z=\bar{q}$, $r=t=q^2$)
 
\n $=\bar{q}([[[a_i^-, a_{i-1}^+],
a_{i+1}^-]_{\bar{q}}, a_{i-1}^+]_{q^2}]_q+\bar{q}[a_{i+1}^-, 
[[ a_i^-, a_{i-1}^+], a_{i-1}^+]_{q^2}]_q)L_{i-1}^2=0$,
according to ($38d$).

\s\n
(iii) $i\in [n+2;n+m].$ Again evaluate first

\n
$
[e_i, e_{i+1}]_q=[\[ a_i^-, a_{i-1}^+\], \[ a_{i+1}^-, 
a_i^+\]]L_iL_{i-1}
$

\n
(from (42) with $a=\[ a_i^-, a_{i-1}^+\]$,  $b=a_{i+1}^-$,
$c=a_i^+$, $x=y=1$, $z=\bar{q}$, $r=s=t=q$) 

\n
$
=\[\[\[ a_i^-, a_{i-1}^+\], a_{i+1}^-\]_{\bar{q}}, a_i^+\]_qL_iL_{i-1}
+\bar{q}\[ a_{i+1}^-, \[\[ a_i^-, a_{i-1}^+\], a_i^+\]_q\]_q L_iL_{i-
1}
$

\n (use ($38d$))

\n
$
=\bar{q}\[ a_{i+1}^-, \bar{L}_ia_{i-1}^+\]
=\[a_{i+1}^-, a_{i-1}^+\]L_{i-1}.
$

\s\n Hence

\n
$
[e_i, [e_i, e_{i+1}]_q]_{\bar{q}}=[\[ a_i^-, a_{i-1}^+\] L_{i-1}, \[ 
a_{i+1}^-,
a_{i-1}^+\]L_{i-1}]_{\bar{q}}
=q[\[a_i^-, a_{i-1}^+\],\[a_{i+1}^-, a_{i-1}^+\]]_{\bar{q}}L_{i-
1}^2
$

\n
(from (42) with $a=\[a_i^-, a_{i-1}^+\]$, $b=a_{i+1}^-$,
$c=a_{i-1}^+$, $x=r=t=1$, $y=z=\bar{q}$, $s=q$ and the triple relations
($38d$))

\n
$=0$. 

The other triple $e$-Serre relation 
$ [e_i, [e_{i}, e_{i-1}]_{\bar{q}}]_q=
  [e_i, [e_{i}, e_{i-1}]_{q}]_{\bar{q}}=0$
is proved in a similar way. 
 
\s\n
F) In order to complete the proof it is convenient
to show as an %%%itermediate 
intermediate step
that Eqs. (28) are consequence of (38). We begin with the 
l.h.s.  of ($28a$). 

\n
$ \[ e_i, a_j^-\]_{q_j^{\delta_{i-1,j}-\delta_{ij}}}=
 \[-(-1)^{\t_{i-1}}\[a_i^-, a_{i-1}^+\]L_{i-1}, 
a_j^-\]_{q_j^{\delta_{i-1,j}-\delta_{ij}}}.
$
Push $L_{i-1}$ to the right and expand the outer supercommutator: 

$$
\eqalign{
& \[ e_i, a_j^-\]_{q_j^{\delta_{i-1,j}-\delta_{ij}}}=
  -(-1)^{\t_{i-1}}\big( q^{1+(-1)^{\t_{i-1}}\delta_{i-1,j}}
\[ a_i^-, a_{i-1}^+\] a_j^-\cr
&-(-1)^{\t_{i,i-1}\t_j}
q_j^{\delta_{i-1,j}-\delta_{ij}}
a_j^-\[ a_i^-, a_{i-1}^+\]\big)L_{i-1}.\cr
}\eqno(44)
$$

\s\n
(i) The case $j<i-1$. From (44)
$
\[ e_i, a_j^-\]
=-(-1)^{\t_{i-1}}q\[\[ a_i^-, a_{i-1}^+\], 
a_j^-\]_{\bar{q}}L_{i-1}=0,
$
according to ($38d$), i.e., ($28a$) holds for $j<i-1$.

\s\n
(ii) The case $j=i-1$. From (44)

\n
$$
\[ e_i, a_{i-1}^-\]_{q_{i-1}}
=-(-1)^{\t_{i-1}}\left( 
q^{1+(-1)^{\t_{i-1}}}\[ a_i^-, a_{i-1}^+\] a_{i-1}^- 
- q_{i-1}a_{i-1}^-\[ a_i^-, a_{i-1}^+\]\right) 
L_{i-1}.
\eqno(45)
$$
(ii.1) If $i\in [1;n+1$], then $\t_{i-1}=0,$ $q_{i-1}=q$ and
(45) \& ($38d$) yield  \hfill\break
$
\[ e_i, a_{i-1}^-\]_q=-q^2\[\[ a_i^-, a_{i-1}^+\] , a_{i-1}^-
\]_{\bar{q}}L_{i-1}=-qa_i^-.
$

\s\n
(ii.2) If $i\in [n+2;n+m]$, then $\t_{i-1}=1,$ $q_{i-1}=\bar{q}$ and
(45) \&  ($38d$)  yield\hfill\break
$
\[ e_i, a_{i-1}^-\]_{\bar{q}}=\[\[a_i^-, a_{i-1}^+\], a_{i-1}^-
\]_{\bar{q}}L_{i-1}=-\bar{q}a_i^-.
$
Hence for $j=i-1$ ($28a$) is fulfilled.

\s\n
(iii) The case $j=i$. Then (44) reduces to
$$
\[ e_i, a_i^-\]_{\bar{q}_i} =-(-1)^{\t_{i-1}}
\left( q\[ a_i^-, a_{i-1}^+\] a_i^- -\bar{q}_i  
a_i^-\[ a_i^-, a_{i-1}^+\]\right) L_{i-1}.\eqno(46)
$$

\n
(iii.1) If $i\in [1;n+1]$, then $\t_{i-1}=0,$ 
$q_{i}=q^{(-1)^{\t_i}}$ and (46)
\& ($38d$) yield \hfill\break
$
\[ e_i, a_i^-\]_{q^{-(-1)^{\t_i}}}=
-q\[\[a_i^-, a_{i-1}^+\], a_i^-\]_{q^{-(1+(-1)^{\t_i})}}L_{i-1}=0.
$

\s\n
(iii.2) If $i\in [n+2;n+m]$, then $\t_{i-1}=1$, $q_{i}=\bar{q}$ 
and (46) \& ($38d$) yield  \hfill\break
$
\[ e_i, a_i^-\]_q=q\[\[a_i^-,a_{i-1}^+\], a_i^-\] L_{i-1}=0. \hfill\break
$
Hence for $i=j$ ($28a$) is fulfilled.

\s\n
(iv) The case $j>i$. Then (44) \& ($38d$) yield \hfill\break
$
\[ e_i, a_{j}^-\]=-(-1)^{\t_{i-1}}q\[\[ a_i^-, a_{i-1}^+\], 
a_j^-\]_{\bar{q}}L_{i-1}=0.\hfill\break
$
Therefore ($28a$) is a consequence of Eqs. (38).

In a similar way one proves that the other relations (28) can be
derived from Eqs. (38).

Note that from Eqs. (28) one derives also Eqs. ($24a,b$).

\s\n
G) We are ready now to derive the additional Serre relation
   $(22c)$.

Using $(24a)$, write $a_{n+2}^-=[[[a_{n-1}^-, e_n]_{\bar{q}}, e_{n+1}
]_{\bar{q}}, e_{n+1}]_q.$ From ($28a$) $\{ e_{n+1}, a_{n+2}^- \}
=0.$ Therefore $0=\{ e_{n+1}, a_{n+2}^- \}=\{ e_{n+1}, [[[a_{n-1}^-, 
e_n]_{\bar{q}}, e_{n+1}]_{\bar{q}}, e_{n+1}]_q\}$. Since $[e_{n+1}, 
a_{n-1}^-]=0,$ and $[e_{n+2}, a_{n-1}^-]=0$ (see ($28a$)) applying 
twice (29) and once (30) one obtains $0=\{ e_{n+1}, a_{n+2}^-\}=
[a_{n-1}^-, y]_{\bar{q}},$ where 

$$
y= \{ e_{n+1},[[e_n,e_{n+1}]_{\bar{q}}, 
e_{n+2}]_{q}\}. \eqno(47)
$$ 
Therefore 
$
[y, a_{n-1}^-]_q=0. 
$
From ($24b$), ($21c$) and  (47)  it follows that $[y, a_{n-
1}^+]=0.$ 
Applying (29) we have 

\noindent
$0=[[y, a_{n-1}^-]_q, a_{n-1}^+]
=[y,[a_{n-1}^-, a_{n-1}^+]]_q$ 

\n(use ($38c$), (24c) and ($21b$)) 

\n$=(q-\bar{q})^{-1}[y, 
L_{n-1}-
\bar{L}_{n-1}]_q=qy\bar{L}_{n-1}.$
Hence,  $y=0$, i.e., the additional $e-$Serre relation ($22c$)
holds.

\s\n
H) In a similar way one derives the $f-$Serre relations $(22d-f)$.
Another way to prove them is to apply the star-operation on
the $e-$ Serre relations.

This completes the proof of the Theorem.

\smallskip
\bigskip\n
{\bf 4. Discussions and further outlook}

\bigskip\n
In the present paper we  enlarge the list of the quantum 
superalgebras, which admit a description via deformed creation and 
annihilation generators [8-13], adding to it all quantum superalgebras
$U_q[sl(n+1|m)]$. The possibility for such a description is not unexpected.
We have generalized the results for $U_q[sl(n+1)]$ [13] to the 
superalgebra case. This generalization is however, we wish to point out,
neither evident nor straightforward. The ``super'' structure is 
richer, with more relations ($e_{n+1}^2=f_{n+1}^2=0$,
additional Serre relations ($22c,f$))  and, as a result, with several features
which do not appear in the Lie algebra cases (the
simple root systems are not related by transformations from
the Weyl group, one and the same superalgebra admits several 
Dynkin diagrams, etc.). All these peculiarities,
especially in the deformed case, which we have mainly in mind here and 
bellow, make the computations nontrivial, technically much more involved.

In the introduction we said few words for a justification of
the name {\it creation and annihilation generators}. Another reason for
this name stems from the observation that, using the CAGs, one
can construct Fock spaces in a much similar way as in the parastatistics
quantum field theory (postulating the existence
of a vacuum, which is annihilated by all $a_i^-$ operators
and introducing an order of the statistics [16];
for more details on parastatistics see, for instance, [32]). Then
the Fock spaces are generated by the creation operators, acting on
the vacuum. Moreover $a_i^+$, acting on a state with fixed
number of ``particles'' (elementary excitation) of species 
$i$, increases them by one, 
whereas $a_i^-$ diminishes them by one. The advantage of this property
for the physical applications and interpretation is evident. 
Consider, for instance, a ``free'' Hamiltonian
$$
H=\sum_{i=1}^{n+m}\varepsilon_i H_i,\quad {\rm such \; that}\quad
  \sum_{i=1}^{n+m}(-1)^{\t_i} \varepsilon_i =0, \eqno(48)
$$
which in the nondeformed case takes the usual form
$$
H=\sum_{i=1}^{n+m}\varepsilon_i \[a_i^+,a_i^-\].\eqno(49)
$$
Then 
$$
[H,a_i^\pm]=\pm \varepsilon_i a_i^\pm,\eqno (50)
$$
i.e., $a_i^+$ (resp. $a_i^-$) can be interpreted as an operator
creating (resp. annihilating) a ``particle'' of species $i$ with
energy $\varepsilon_i$.
Our, we call it {\it physical conjecture} is that the Fock
representations of the deformed CAGs will lead to new solutions 
for the microscopic  $g-$ons statistics in the sense of Karabali and 
Nair [33], which is a particular realization of the exclusion statistics
of Haldane [27].

The Fock representations 
however may be of interest also from another, more mathematical point 
of view. So far the finite-dimensional irreducible representations of the
LSs from the class $A$ were explicitly constructed only for 
$sl(n|1)$ [34]. Any such representation 
can be deformed to a representation of $U_q[sl(n|1)]$ [35].
The representation theory of $sl(n|m),\;n,m=1,2,\l$ and hence of the
corresponding deformed algebras is however far from being
complete, if both  $n\ne 1$ and $ m\ne 1$. In [36] the so-called
essentially typical representations of $sl(n|m)$ were
described. The results were generalized also to the quantum 
case [37]. Our {\it mathematical conjecture} now is that the Fock 
representations are beyond the class of the deformed 
essentially typical representation [36], thus yielding new 
representations of $U_q[sl(n+1|m)]$.

In order to verify the above conjectures one would need
to construct the Fock representations explicitly, i.e., to introduce a 
basis and to write down the transformations of the basis under the
action of the generators. As a first step one has to determine 
the quantum analogue of the 
triple relations (17). This is a nontrivial problem.
It actually means that one has to write down the supercommutation
relations between all Cartan-Weyl generators, expressed via
the CAGs. The latter is a necessary condition for the application of
the Poincare-Birghoff-Witt theorem, when computing the action of
the generators on the Fock basis vectors. We return to this problem 
elsewhere. Here we mention only one, but important additional relation: 
from (17) one derives that the creation (resp. annihilation) 
generators $q$-supercommute,
$$
 \[ a_i^\xi , a_j^\xi \]_{q'} =0,\quad q'=q\;{\rm or}\;\q,\;
\; i,j\in [1;n+m], \;\; \xi=\pm.\eqno(51)
$$
This makes evident the basis (or at least one possible basis) in a 
given Fock space, since any product of only creation generators 
can be always ordered. Note that similar property does not hold for
para-Bose (or para-Fermi) creation operators. This is  
the reason why (even in the nondeformed case) 
the matrix elements of the paraoperators
remain still unknown for an arbitrary order of the
parastatistics: the Fock space basis cannot be represented as ordered
products of only para-Bose (or para-Fermi) creation 
operators acting on the vacuum (the linear span of  only
such vectors is not invariant under the action of the para-operators). 

Finally let us mention that we do not have simple
relations for the action of the other Hopf algebra
operations $(\Delta, \varepsilon, S)$ on the CAGs,
although it is clear how to write them down, using Eqs. (16)
and the circumstance that the comultiplication $\Delta$ and
the counity $\varepsilon$ are morphisms, whereas the antipode
$S$ is an antimorphism. In this respect the picture is much 
the same as discussed in [13]. Luckily, the 
$(\Delta, \varepsilon, S)$-operations are not necessary for
%%%computig
computing  the transformations of the Fock modules (but they are
certainly very important when considering tensor products
of representation spaces).

\vskip 20pt
\noindent
{\bf Acknowledgments.}

\s\n
We are thankful to Prof. Randjbar-Daemi for the kind invitation
to visit the High Energy Section of the Abdus Salam International 
Centre for Theoretical Physics. 

\b

\n {\bf References}

\s\n
{\settabs\+[11] & I. Patera, T. D. Palev, Theoretical interpretation of the 
   experiments on the elastic \cr 
   %sample line,  see p. 232 of the Texbook.

\+ [1] & Drinfeld V G 1985 {\it DAN SSSR} {\bf 283} 1060;
         1985 {\it Sov. Math. Dokl.} {\bf 32} 254 \cr
\+ [2] & Jimbo M 1985 {\it Lett. Math. Phys.}{\bf 10} 63 \cr 

\+ [3] & Kulish P P 1988 {\it Zapiski nauch. semin. LOMI}
         {\bf 169} 95 \cr
\+ [4] & Kulish P P and  Reshetikhin N Yu 1989 {\it Lett. Math.
         Phys.} {\bf 18} 143 \cr
\+ [5] & Chaichian M and Kulish P P 1990 {\it Phys. Lett.} 
         {\bf B 234} 72\cr
\+ [6] & Bracken A J, Gould M D and Zhang R B 1990 
         {\it Mod. Phys. Lett.} {\bf A 5} 831 \cr
\+ [7] & Tolstoy V N 1990 {\it Lect. Notes in Physics} {\bf 370},
         Berlin, Heidelberg, New York:\cr
\+     & Springer p. 118 \cr
\+ [8] & Palev T D 1993 {\it J. Phys. A: Math. Gen.}
         {\bf 26} L1111 and hep-th/9306016 \cr
\+ [9] & Hadjiivanov L K 1993 {\it J. Math. Phys} {\bf 34} 5476 \cr
         
\+ [10] & Palev T D and Van der Jeugt 1995 {\it J. Phys. A: Math. Gen.}
         {\bf 28} 2605 \cr
\+      & and q-alg/9501020 \cr
\+ [11] & Palev T D 1984 {\it Lett. Math. Phys.} {\bf 31} 
          151  and het-th/9311163 \cr
\+ [12] & Palev T D 1998 {\it Commun. Math. Phys.} {\bf 196} 429
          and q-alg/9709003 \cr
\+ [13] & Palev T D and Parashar P 1998 {\it Lett. Math. Phys.}
          {\bf 43} 7 and q-alg/9608024   \cr
\+ [14] & Palev T D 1980 {\it J. Math. Phys} {\bf 21} 1293 \cr

\+ [15] & Palev T D 1982 {\it J. Math. Phys} {\bf 23} 1100 \cr

\+ [16] & Green H S 1953 {\it Phys. Rev.} {\bf 90} 270\cr 
\+ [17] & Kac V G 1979 {\it Lecture Notes in Math.} {\bf 676},
           Berlin, Heidelberg, New York:\cr
\+      &  Springer p. 597 \cr
\+ [18] & Palev T D 1979 {\it Czech. Journ. Phys.} {\bf B 29} 91\cr
\+ [19] & Palev T D 1976 Lie algebraical aspects of the quantum
          statistics {\it Thesis}  \cr
\+      & Institute for Nuclear Research and Nuclear Energy, 
          Sofia \cr          
\+ [20] & Palev T D 1977 Lie algebraical aspects of the quantum
          statistics. Unitary \cr
\+      & quantization (A-quantization) 
          {\it Preprint JINR E17-10550} and hep-th/9705032;\cr
\+      & 1980 {\it Rep. Math. Phys} {\bf 18} 117 and 129   \cr
\+ [21] & Palev T D 1978 A-superquantization {\it Communications JINR}
          E2-11942\cr
\+ [22] & Palev T D 
          1982 {\it J. Math. Phys.} {\bf 23} 1778;
          1982 {\it Czech. Journ. Phys.} {\bf B 32} 680 \cr

\+ [23] & Okubo S 1994 {\it J. Math. Phys.} {\bf 35} 2785  \cr

\+ [24] & Van der Jeugt J 1996  {\it New Trends in Quantum Field
          Theory} (Heron Press, Sofia) \cr
\+ [25] & Meljanac S, Milekovic M and Stojic M 1998 {\it Mod. Phys. Lett}
          {\bf A 13} 995 \cr
\+      & and q-alg/9712017 \cr
\+ [26] & Palev T D and Stoilova N I 1997 {\it J. Math. Phys.} 
          {\bf 38} 2506 and hep-th/9606011 \cr
\+ [27] & Haldane F D M 1991 {\it Phys. Rev. Lett.} {\bf 67} 937\cr
\+ [28] & Palev T D 1992 {\it Rep. Math. Phys} {\bf 31} 241 \cr
\+ [29] & Khoroshkin S M and Tolstoy V N 1991 {\it Commun. Math. Phys.} 
          {\bf 141} 599  \cr 
\+ [30] & Floreanini R Leites D A and Vinet L 1991 {\it Lett. Math. Phys.}
          {\bf} 23 127 \cr
\+ [31] & Scheunert M 1992 {\it Lett. Math. Phys.} {\bf 24} 173 \cr 

\+ [32] & Ohnuki Y and Kamefuchi S 1982 {\it Quantum Field Theory
          and parastatistics} \cr
\+      & (University of Tokyo, Springer Verlag)\cr

\+ [33] & %%%Karballi
         Karabali D and Nair V P 1995 {\it Nucl. Phys.}
          {\bf B 438} 551\cr

\+ [34] & Palev T D 1989 {\it J. Math. Phys.} {\bf 30} 1433 \cr

\+ [35] & Palev T D and Tolstoy V N 1991 {\it Commun. Math. Phys.}
          {\bf 141} 549 \cr

\+ [36] & Palev T D 1989 {\it Funkt. Anal. Prilozh.} {\bf 23}
          $\# 2$ 69 (in Russian);\cr

\+      & 1989 {\it Funct. Anal. Appl.}
          {\bf 23} 141 (English translation)\cr

\+ [37] & Palev T D Stoilova N I and Van der Jeugt J 1994
          {\it Commun. Math. Phys.} {\bf 166} 367 \cr

\end